# An Empirical Investigation of
# Four Well-Known Polynomial-Size VRP Formulations

Deniz Aksen [a,*], Temel Öncan [b], Mir Ehsan Hesam Sadati [c]

[a] Koç University; College of Administrative Sciences and Economics
Rumelifeneri Yolu, Sarıyer, 34450 İstanbul, Türkiye
daksen@ku.edu.tr

[b] Galatasaray University; Department of Industrial Engineering
Çırağan Cad. No:36, Ortaköy, 34357 İstanbul, Türkiye
ytoncan@gsu.edu.tr

[c] Koç University; Department of Industrial Engineering
Rumelifeneri Yolu, Sarıyer, 34450 İstanbul, Türkiye
msadati14@ku.edu.tr

**Abstract:**

This study presents an in-depth computational analysis of four well-known Capacitated Vehicle Routing Problem (CVRP) formulations with polynomial number of subtour elimination constraints: a node-based formulation and three arc-based (single, two- and multi-commodity flow) formulations. For each formulation, several valid inequalities (VIs) are added for the purpose of tightening the formulation. Moreover, a simple topology-driven granulation scheme is proposed to reduce the number of a certain type of VIs. The lower and upper bounding performance and the solution efficiency of the formulations and respective VI configurations are benchmarked with state-of-the-art commercial optimization software. The extensive computational analysis embraces 121 instances with up to 100 customer nodes. We believe that our findings could be useful for practitioners as well as researchers developing algorithms for the CVRP.

## 1. Introduction

During the last six decades, the Capacitated Vehicle Routing Problem (CVRP) has attracted the interest of many researchers, and numerous papers have been published ever since the seminal paper by Dantzig and Ramser (1959). A simple phrase match search for "vehicle routing" covering the years from 2006 through 2016 returns about 24,600 results in Google Scholar[1]. The CVRP and its extensions have a substantial number of real-world applications in logistics (Toth and Vigo, 2002). The CVRP is known to be $\mathcal{NP}$-hard (Lenstra and Rinnooy Kan, 1981) and can be represented as a graph theoretic problem. Let $G = (I, A)$ be a complete directed graph, where $I = \{1, ..., n\}$ is the node set and $A$ is the arc set. Each node $i = 1, ..., n$ corresponds to a customer and node 0 stands for the depot. A nonnegative cost $d_{ij}$ corresponds for each arc $(i, j) \in A$, which implies the travel distance from node $i$ to node $j$. Note that, in this study we address the asymmetric case of the CVRP where $d_{ij} \neq d_{ji}$ holds for $i, j = 0, ..., n; i \neq j$ and we assume $d_{ii} = \infty$ for $i = 0, ..., n$. The CVRP consists of finding $K$ vehicle tours, i.e., $K$ cycles on $G$, with the minimum total travel distance starting from depot (node 0). The problem is solved for a fleet of homogenous vehicles with a fixed capacity $Q$ which visit $n$ customers with nonnegative demand quantities $q_i$ for $i=1,…, n$. Furthermore, we assume that each customer is visited by exactly one vehicle, the load of the vehicle should not exceed the vehicle capacity $Q$ and each route starts and ends at the depot node 0.

---

[*] Corresponding author. Tel.: +90 (1) 338 1684, Fax.: +90 (1) 338 1653, Email: daksen@ku.edu.tr (D. Aksen)

[1] Google Scholar search results obtained as of Aug 20, 2016 for the query enclosed in double quotation marks.

We intend to accomplish an in-depth computational analysis and comparison of four CVRP formulations with polynomial number of subtour elimination constraints: a node-based formulation using load variables (Desrochers and Laporte, 1991) and three arc-based formulations using flow variables (Gavish and Graves, 1978; Baldacci, Mingozzi and Hadjiconstantinou, 2004, Letchford and Salazar-González, 2006 & 2015). The strengths of node-based and arc-based formulations have been discussed by several researchers in the context of the Traveling Salesman Problem (TSP). Examples include Padberg and Sung (1991), Öncan, Altınel and Laporte (2009), Roberti and Toth (2012). Recently, Sarin et al. (2014) investigate the performances of 32 formulations for the Multiple Asymmetric Traveling Salesman Problem with and without precedence constraints (mATSP and PCmATSP). The authors compare the tightness of the formulations, and assess their solvability using commercial software. The results suggest that the tightest LP relaxations are afforded by the arc-based formulations. In the realm of the CVRP we are aware of only one paper by Ordónez, Sungur and Dessouky (2007). They conduct experiments addressing the effect of several problem characteristics such as the capacity $Q$, the number of vehicles $K$ and the distance matrix on the performance of the unit-demand CVRP. To our knowledge, there is no other study comparing the effectiveness and efficiency of both node-based and arc-based formulations for the CVRP.

Our empirical findings are based on extensive and systematic experimentation. We believe that they can be projected onto the modeling of a wide range of optimization problems involving routing decisions. The performance metrics which have been meticulously evaluated can be utilized as a benchmark criterion for many solution approaches. This research may appeal also to practitioners seeking quick and easy-to-obtain solutions with small, yet proven optimality gaps.

Our motivation is to demonstrate that—using modern-day computer hardware and state-of-the-art optimization software—the right combination of a mathematical model, a set of valid inequalities, and multithreading options of the right solver can make a sizable difference. Such a combination can obtain favorable solutions to CVRP instances with up to 100 customers within reasonable times. However, we should state that a rigorous theoretical comparison of CVRP formulations is beyond the scope of this paper.

The remainder of the paper is organized as follows. In Section 2 we present four polynomial-size CVRP formulations. The experimental setup and the results obtained in the experiments are discussed in Section 3. Finally, conclusions and future research are made in Section 4.

## 2. Node-based and arc-based formulations for the CVRP

The node-based formulation considered in this study was originally suggested by Desrochers and Laporte (1991) for the TSP. It uses a lifted version of the well-known Miller-Tucker-Zemlin (MTZ) inequalities that had been first proposed by Miller, Tucker and Zemlin (1960) for the elimination of subtours in TSP solutions. In addition to this node-based formulation, we have also considered three arc-based formulations: the first one is the single-commodity network flow formulation developed by Gavish and Graves (1978) for the TSP. The second one is the two-commodity network flow formulation developed for the CVRP by Baldacci, Mingozzi and Hadjiconstantinou (2004). Finally, the third one is the multi-commodity flow (MCF) formulation of CVRP which is proposed by Letchford and Salazar-González (2006, 2015).

For the TSP case, although the Desrochers and Laporte (1991) formulation is incomparable to the Gavish and Graves (1978) formulation (Öncan, Altınel and Laporte, 2009), Roberti and Toth (2012) have reported that the model devised by Desrochers and Laporte (1991) arises to be computationally competitive among several polynomial-size formulations of the asymmetric TSP. Hence, we believe it is an open research question to compare the computational performance of node-based and arc-based formulations for the CVRP. We remark that the CVRP has been modeled in the first two formulations coming next as a collection problem rather that a delivery problem. Changing this convention for both node-based and arc-based formulations is straightforward, and is not expected to make a difference in the solution speed or accuracy of the model in question.

Let us define the index set $\mathbf{I}_0 = \{0, 1, \ldots, n\}$ as the set of $n$ customer nodes and the depot 0, and $\mathbf{I} = \{1, \ldots, n\}$ as the set of $n$ customer nodes only. Furthermore, the following parameters are defined. Let $d_{ij}$ be the distance from node $i$ to node $j$, $(i, j \in \mathbf{I}_0, d_{ij} \neq d_{ji})$, $Q$ be the vehicle capacity in demand units and, $q_i$ be the demand or load to be delivered to or to be collected from customer node $i$, $(i \in \mathbf{I})$. The decision variable used in both node-based and arc-based formulations is as follows. Let $X_{ij}$ be a binary variable indicating if arc $(i, j)$ is traversed by a vehicle $(i, j \in \mathbf{I}_0)$. Finally let $Z$ be the total traveling distance or cost (objective function value). A generic CVRP formulation can be given as follows:

$$\text{minimize } Z = \sum_{i \in \mathbf{I}_0} \sum_{j \in \mathbf{I}_0,\, i \neq j} d_{ij} X_{ij} \qquad (1)$$

subject to

$$\sum_{j \in \mathbf{I}_0,\, i \neq j} X_{ji} = 1 \qquad i \in \mathbf{I} \qquad (2)$$

$$\sum_{j \in \mathbf{I}_0, i \neq j} X_{ij} = 1 \qquad i \in \mathbf{I} \qquad (3)$$

$$\text{Subtour Elimination } + \text{ Feasible Tour Construction Constraints} \qquad (4)$$

$$X_{ij} \in \{0,1\} \qquad i, j \in \mathbf{I}_0 \qquad (5)$$

where assignment-like constraints (2) and (3) state that each customer has exactly one incoming and one outgoing vehicle. Constraints (4) state that each vehicle performs a feasible tour, i.e. the vehicle capacity is not exceeded, and there are no subtours. Finally, constraints (5) guarantee the integrality of the decision variables.

## 2.1. A node-based formulation for the CVRP

Our first formulation uses node-specific decision variables $U_i$ ($i \in \mathbf{I}$) which indicate the load of a vehicle right after departing from customer $i$. We write the following Lifted MTZ inequalities for subtour elimination in the CVRP modeled as a collection problem.

$$U_i - U_j + Q X_{ij} + (Q - q_i - q_j) X_{ji} \leq Q - q_j \qquad i, j \in \mathbf{I},\ i \neq j \qquad (6)$$

$$U_i \geq q_i + \sum_{j \in \mathbf{I}, j \neq i} q_j X_{ji} \qquad i \in \mathbf{I} \qquad (7)$$

$$U_i \leq Q - \sum_{j \in \mathbf{I}, j \neq i} q_j X_{ij} \qquad i \in \mathbf{I} \qquad (8)$$

$$U_i \leq Q - (Q - q_i) X_{0i} \qquad i \in \mathbf{I} \qquad (9)$$

$$U_i \leq Q - (Q - \max_{\substack{j \in \mathbf{I} \\ j \neq i}}\{q_j\} - q_i) X_{0i} - \sum_{j \in \mathbf{I}, j \neq i} q_j X_{ij} \qquad i \in \mathbf{I} \qquad (10)$$

$$U_i \geq 0 \qquad i \in \mathbf{I} \qquad (11)$$

The validity of constraints (6)–(8) follows from the definition of the variables $U_i$. The validity of constraints (6) can be demonstrated as follows: Since $U_i$ stands for the load of a vehicle right after departing from customer $i$ then $U_i$ should not be less than the load carried with the vehicle coming from customer $j$ plus the load collected from customer $i$. On the other hand, due to capacity constraints, the values of $u_i$ variables cannot exceed the vehicle capacity $Q$ unless the vehicle travels from customer $i$ to customer $j$. In that case, i.e. when $X_{ij} = 1$, then $U_i \leq Q - q_j X_{ij}$ holds which implies the validity of constraints (6). Similarly, one can also show that constraints (9) are valid for the CVRP. By definition $U_i \leq Q$ holds for all customers $i \in \mathbf{I}$. However, when $X_{0i} = 1$ holds, then the vehicle immediately travels from depot 0 to customer $i$, hence by definition we have $U_i = q_i$ which shows the validity of the constraints (9). Finally, to demonstrate the validity of the constraints (10) we will again use the definition of $U_i$. In case there is no incoming vehicle from depot to customer $i$ and there is no outgoing vehicle from customer $i$ to customer $s$, then $U_i \leq Q$ holds. In case, the vehicle travels from depot to customer $i$ and from customer $i$ to customer $j$ then $U_i \leq Q - (Q - \max_{\substack{j \in \mathbf{I} \\ j \neq i}}\{q_j\} - q_i) - q_s = \max_{\substack{j \in \mathbf{I} \\ j \neq i}}\{q_j\} + q_i - q_s$ holds. Besides, $U_i \leq \max_{\substack{j \in \mathbf{I} \\ j \neq i}}\{q_j\} + q_i - q_s$ is always satisfied even if $\max_{\substack{j \in \mathbf{I} \\ j \neq i}}\{q_j\} = q_s$ holds by the definition.

We will use the terms node-based formulation and MTZ-L interchangeably to denote the model that consists of the equations (1)–(3), integrality constraints (5), and the Lifted MTZ inequalities in (6)–(11).

## 2.2. Single-commodity flow formulation for the CVRP

Our second model is derived from the single-commodity flow formulation proposed initially for the TSP in Gavish and Graves (1978). It uses arc-specific decision variables $F_{ij}$ ($i \in \mathbf{I}$, $j \in \mathbf{I}$) to indicate the amount of goods flowing from customer $i$ to customer $j$ after collecting the load at customer $i$. In order to assure subtour elimination in the formulation, we impose tight bounds on the flow variables $F_{ij}$ in addition to a set of flow conservation constraints as follows.

Flow conservation at each customer node: $\sum_{j \in \mathbf{I}_0, j \neq i} F_{ji} + q_i = \sum_{j \in \mathbf{I}_0, j \neq i} F_{ij}$ $\quad i \in \mathbf{I}$ (12)

Lower bound on the $F_{ij}$: $\quad F_{ij} \geq q_i X_{ij}$ $\quad i, j \in \mathbf{I}_0, i \neq j$ (13)

Upper bound on the $F_{ij}$: $\quad F_{ij} \leq (Q - q_j) X_{ij}$ $\quad i, j \in \mathbf{I}_0, i \neq j$ (14)

Integrality and non-negativity constraints: $\quad F_{ij} \geq 0$ $\quad i, j \in \mathbf{I}_0$ (15)

The single-commodity flow model that consists of equations (1)–(3), (5), (12)–(15) will be denoted as GG in recognition of the introduction of this type of formulations in Gavish and Graves (1978).

## 2.3. Two-commodity flow formulation for the symmetric CVRP

Our third CVRP model is due to Baldacci, Hadjiconstantinou and Mingozzi (2004). It is a successful adaptation of the two-commodity flow formulation of the TSP which was first proposed by Finke, Claus and Gunn (1984). The two-commodity flow formulation works basically for the symmetric CVRP modeled as a delivery problem with a fixed number of vehicles. Let us define the index set $\mathbf{I}'_0 = \{0, 1, \ldots, n, n+1\}$ obtained from $\mathbf{I}_0$ by adding a dummy node $(n+1)$ where $(n+1)$ is an identical copy of the depot. Let us define now two nonnegative flow variables $G_{ij}$ and $G_{ji}$ to represent an edge $(i, j)$ of a feasible CVRP solution on the original graph augmented by the dummy depot. The flow variables $G_{ij}$ and the routing variables $X_{ij}$ are defined for all edges but $(0, n+1)$. The vehicle is assumed to carry a combined load of $Q$ units. For $i$ and $j \in \mathbf{I}'_0$ ($i \neq j$), the flow variables $G_{ij}$ specify two directed flow paths in opposite directions for any route of a feasible solution. One path goes from node 0 to node $(n + 1)$; the vehicle assigned to this path is dispatched from the depot node 0 with the total demand of customer nodes to be visited on it. Another path goes in the opposite direction from $(n + 1)$ to 0 with the total flow equal to the vehicle capacity $Q$. If $X_{ij}$ is equal to 1, i.e. if a vehicle travels from node $i$ to node $j$, then $G_{ij}$ represents the load of the vehicle by the time it departs from node $i$, and $G_{ji}$ indicates the residual capacity of the vehicle on the same edge $(i, j)$. Note that $G_{ji} = Q - G_{ij}$ holds true whenever $X_{ij} = 1$. The two-commodity flow formulation for the CVRP is given as follows in Baldacci, Hadjiconstantinou and Mingozzi (2004):

$$\text{minimize } Z' = \sum_{i \in \mathbf{I}_0} \sum_{j \in \mathbf{I}'_0, i \neq j} d_{ij} X_{ij} + \sum_{i \in \mathbf{I}} d_{i(n+1)} X_{i(n+1)} \tag{16}$$

subject to

$$\sum_{j \in \mathbf{I}'_0} \left( G_{ji} - G_{ij} \right) = 2q_i \qquad \forall i \in \mathbf{I} \tag{17}$$

$$\sum_{j \in \mathbf{I}} G_{0j} = \sum_{j \in \mathbf{I}} q_j \tag{18}$$

$$\sum_{j \in \mathbf{I}} G_{j0} = KQ - \sum_{j \in \mathbf{I}} q_j \tag{19}$$

$$\sum_{j \in \mathbf{I}} G_{(n+1)j} = KQ \tag{20}$$

$$G_{ij} + G_{ji} = QX_{ij} \qquad \forall i, j \in \mathbf{I}'_0, \ i \neq j, i \neq (n+1) \tag{21}$$

$$\sum_{\substack{j \in \mathbf{I}'_0 \\ j > i}} X_{ij} + \sum_{\substack{j \in \mathbf{I}'_0 \\ j < i}} X_{ji} = 2 \qquad \forall i \in \mathbf{I} \tag{22}$$

$$G_{ij} \geq 0 \qquad \forall i, j \in \mathbf{I}'_0 \tag{23}$$

$$X_{ij} \in \{0,1\} \qquad \forall i, j \in \mathbf{I}'_0 \tag{24}$$

The objective function $Z'$ in (16) is exactly twice the one of the previous models in (1). Constraints (17) are for the flow conservation of the two commodities entering and leaving customer node $i$. Constraints (18)$-$(20) make sure that the required amount of commodities enter and leave the depot. Constraint (18) states that the outflow at the depot node 0 should be equal to the total customer demand, and constraint (20) implies that the inflow at the dummy depot node $(n + 1)$ should be equal to the total capacity of the fixed-size vehicle fleet. Consequently, constraint (19) sets the inflow at the depot node 0 equal to the residual capacity of the vehicle fleet. Constraints (21) ensure that when $X_{ij} = 1$ then $G_{ji} = Q - G_{ij}$ holds. Constraints (22) guarantee that there exist exactly two arcs incident to each customer node $i$. Finally, constraints (23) and (24) indicate the non-negativity and binary constraints on the commodity flow variables and the routing variables, respectively.

We denote the two-commodity flow formulation consisting of the equations (16)–(24) as BHM henceforth in recognition of Baldacci, Hadjiconstantinou and Mingozzi's paper (2004).

## 2.4. Multi-commodity flow formulation for the CVRP

Our fourth model belongs to the class of multi-commodity flow formulations. We take into account a particular formulation called "MCF2a" which has a polynomial number of subtour elimination constraints (SECs). MCF2a was first proposed in Letchford and Salazar-González (2006). It uses two extra decision variables $F_{ij}^k$ and $G_{ij}^k$ (representing two commodities per customer). $F_{ij}^k$ is a binary variable defined for each customer $k \in \mathbf{I}$ and each arc $(i, j) \in A$ taking the value 1 if and only if a vehicle traverses $(i, j)$ on the way from the depot to $k$. Likewise, $G_{ij}^k$ is another binary variable defined for each customer $k \in \mathbf{I}$ and each arc $(i, j) \in A$ taking the value 1 if and only if a vehicle traverses $(i, j)$ on the way from $k$ to the depot. MCF2a treats the depot either as the source or as the sink of every commodity. The same authors Letchford and Salazar-González (2015) recently modified this formulation further to obtain what they name as "MCF2b". They prove that MCF2b has the strongest LP relaxation, thus yields the tightest lower bounds among all known flow formulations of polynomial size. Their computational results support this proposition as well. We adopt MCF2b as our multi-commodity flow formulation of choice, and designate it simply as "MCF" in the sequel. According to Letchford and Salazar-González (2015), MCF can be written by replacing (4) with the following set of constraints:

$$\sum_{j \in \mathbf{I}} F_{0j}^k - \sum_{j \in \mathbf{I}} F_{j0}^k = 1 \qquad k \in \mathbf{I} \qquad (25)$$

$$\sum_{j \in \mathbf{I_0}, j \neq i} F_{ij}^k - \sum_{j \in \mathbf{I_0}, j \neq i} F_{ji}^k = 0 \qquad i, k \in \mathbf{I} \quad i \neq k \qquad (26)$$

$$\sum_{j \in \mathbf{I}} G_{j0}^k - \sum_{j \in \mathbf{I}} G_{0j}^k = 1 \qquad k \in \mathbf{I} \qquad (27)$$

$$\sum_{j \in \mathbf{I_0}, j \neq i} G_{ij}^k - \sum_{j \in \mathbf{I_0}, j \neq i} G_{ji}^k = 0 \qquad i, k \in \mathbf{I} \quad i \neq k \qquad (28)$$

$$F_{ij}^k + G_{ij}^k \leq X_{ij} \qquad k \in \mathbf{I}, (i, j \in \mathbf{I_0}) \quad i \neq j \qquad (29)$$

$$\sum_{k \in \mathbf{I} \setminus \{i.j\}} q_k (F_{ij}^k + G_{ij}^k) \leq (Q - q_i - q_j) X_{ij} \qquad (i, j \in \mathbf{I_0}) \quad i \neq j \qquad (30)$$

$$F_{ij}^k, G_{ij}^k \geq 0 \qquad k \in \mathbf{I}, (i, j \in \mathbf{I_0}) \qquad (31)$$

## 2.5. Several straightforward valid inequalities for the CVRP

In this section we discuss several straightforward and well-known valid inequalities that can be applied to the CVRP. The reader is referred to Toth and Vigo (2002) for a more elaborate discussion on this topic.

### 2.5.1 Bounds on the number of arcs leaving the depot

The minimum number of vehicles constraint designated as MinNV in (32) and the maximum number of vehicles constraint designated as MaxNV in (33) below are applicable to the MTZ-L, GG as well as MCF formulations of the CVRP. Let us consider the following consolidated total vehicle capacity constraint which can be treated as a trivial valid inequality:

**MinNV**: $\quad Q \sum_{i \in \mathbf{I}} X_{0i} \geq \sum_{i \in \mathbf{I}} q_i$ (32)

Constraint (32) stipulates that the total demand is satisfied by dispatching from the depot a sufficient number of vehicles of identical capacity $Q$. In order to see that (32) is valid, suppose that $\sum_{i \in \mathbf{I}} X_{0i} < \frac{1}{Q} \sum_{i \in \mathbf{I}} q_i$ holds. Then for at least one tour the vehicle capacity will be violated which implies the validity of constraint (32). We designate constraint (32) as the minimum number of vehicles constraint MinNV in the sequel. Note that we experimented also with a tighter version of constraint (32) which requires $\sum_{i \in \mathbf{I}} X_{0i} \geq \left\lceil \sum_{i \in \mathbf{I}} q_i / Q \right\rceil$ hold true. We compared the CPU times of the GG models with tight and loose versions of MinNV on 15 test instances in which the problem size $n$ varied between 41 and 50 customers, and the loose-MinNV version of the respective GG model was solved to proven optimality within three hours. A two-tailed paired t-test with $\alpha = 0.05$ revealed that there was no significant difference in the CPU times between loose and tight MinNV constraints, $t(14) = -1.5405, \ p = 0.1457$.

The counterpart of MinNV, namely MaxNV sets an upper bound on the number of outgoing arcs at the depot. MaxNV can be written only for those instances where the problem data specifies the size of the available vehicle fleet. In those instances the number of vehicles is set to a fixed $K$ value. We solve 104 such instances as can be seen in the repository in Table A.1. MaxNV is formulated as follows:

**MaxNV**: $\quad \sum_{i \in \mathbf{I}} X_{0i} \leq K$ (33)

The two-commodity flow formulation BHM presented in Section 2.3 can also be supplemented with a similar constraint. Recall that the BHM model works only in the case of a fixed-size fleet of identical vehicles. This feature can be enforced directly by an additional equality constraint as follows.

**FixedK**: $\quad \sum_{i \in \mathbf{I}} X_{0i} = K$ (34)

In order to show the validity of the constraint (34), we write the constraints (21) of the BHM model for $i = 0$ and sum them up over all customers $j \in \mathbf{I}$ which yields the equality $\sum_{j \in \mathbf{I}} G_{0j} + \sum_{j \in \mathbf{I}} G_{j0} = Q \sum_{j \in \mathbf{I}} X_{0j}$. The left-hand side of this equality is the same as the sum of the right-hand sides of constraints (18) and (19), namely $KQ$. We substitute $KQ$ in the left-hand side of the equality, and then divide both sides by the nonzero capacity value $Q$. This simplification leads to the FixedK constraint in (34).

### 2.5.2 Three valid inequalities applicable to MTZ-L and GG

Let us now look into three specific valid inequalities (VIs) marked as VI-1 through VI-3. These VIs can be embedded into the formulations MTZ-L and GG, but cannot be applied to the BHM formulation. We first take into account the following single equation which balances the number of leaving and entering arcs at the depot (node 0).

**Valid Inequality VI-1**:
$$\sum_{i \in \mathbf{I}} X_{0i} = \sum_{i \in \mathbf{I}} X_{i0} \tag{35}$$

The validity of constraint (35) follows from the MinNV constraint in (32) and the assignment-like constraints in (2) and (3). If there are $\sum_{i \in \mathbf{I}} X_{0i}$ vehicles leaving the depot, then due to the assignment constraints there must be exactly the same number of vehicles returning to the depot.

Recall that the subtour elimination constraints originally devised by Dantzig, Fulkerson and Johnson (1954) for the TSP are as follows:

$$\sum_{(i,j) \in S} X_{ij} \leq |S| - 1 \quad \text{for} \quad 2 \leq |S| < n - 1, \; S \subseteq V \tag{36}$$

In the context of CVRP, these inequalities can be extended into the Generalized Subtour Elimination Constraints (GSECs):

$$\sum_{i \in S} \sum_{\substack{j \in S \\ j \neq i}} X_{ij} \leq |S| - L(S) \quad \text{for} \; S \subseteq IC; S \neq \emptyset \tag{37}$$

where is $L(S)$ stands for the minimum number of vehicles needed to enter or leave the customer subset $S$. $L(S)$ can be replaced with a trivial Bin Packing Problem lower bound such as $L(S) = \left\lceil \frac{q(S)}{Q} \right\rceil$ where $q(S)$ is the total demand of customers in subset $S$. It is straightforward to derive valid inequalities $X_{ij} + X_{ji} \leq 2 - \left\lceil \frac{q_i + q_j}{Q} \right\rceil$ for $S = \{i,j\}$ and $X_{ij} + X_{ji} + X_{ik} + X_{ki} + X_{kj} + X_{jk} \leq 3 - \left\lceil \frac{q_i + q_j + q_k}{Q} \right\rceil$ for $S = \{i,j,k\}$. So, we write the following GSECs of size two and size three designated as VI-2 and VI-3, respectively, for the symmetric CVRP:

**Valid Inequality VI-2**:
$$X_{ij} + X_{ji} \leq 2 - \left\lceil \frac{q_i + q_j}{Q} \right\rceil \qquad i, j \in \mathbf{I}, \; i < j \tag{38}$$

**Valid Inequality VI-3**: $\quad X_{ij} + X_{jk} + X_{ki} \leq 3 - \left\lceil \dfrac{q_i + q_j + q_k}{Q} \right\rceil \quad i,j,k \in \mathbf{I},\ i < j < k \quad (39)$

### 2.5.3 Valid equalities applicable to MCF

The LP relaxation of the MCF formulation satisfies the following equations designated as FGX, which indicate the coupling between the variables **F**, **G** and **X**:

$$\mathbf{FGX:} \quad F_{ik}^{k} = X_{ik} \quad i,k \in \mathbf{I}\ i \neq k$$

$$G_{ik}^{i} = X_{ik} \quad i,k \in \mathbf{I}\ i \neq k \quad (40)$$

For the proof of the validity of FGX the reader is referred to Lemma 2 in Letchford and Salazar-González (2015).

## 2.6. Handling the maximum tour length or tour duration

In our experiments we solve four CVRP instances with a maximum tour duration constraint (*MTD*). This constraint requires that the total travel time, i.e. the time between the departure and return of each vehicle be restricted by a maximum tour duration $D_{\max}$. Usually there exists a fixed or customer-dependent service time $fst_i$ at each node $i \in \mathbf{I}$. Travel times are indicated by the entries of the matrix $[t_{ij}]$. If $fst_i$ values are ignored and travel times $t_{ij}$ are replaced by travel distances $d_{ij}$, *MTD* becomes the maximum tour length constraint. We explain the additional decision variables and equations needed to incorporate this constraint into the CVRP model. *MTD* is handled in the same way in both the node-based and the arc-based formulations.

Let the continuous decision variable $A_i$ denote the arrival time of some vehicle at the customer $i$. We write two sets of equations to impose lower and upper bounds on $A_i$ and one set of equations to impose the *MTD*.

$$A_j \geq A_i + fst_i + t_{ij} - (1 - X_{ij})M \quad i \in \mathbf{I_0},\ j \in \mathbf{I},\ i \neq j \quad (41)$$

$$A_j \leq A_i + fst_i + t_{ij} + (1 - X_{ij})M \quad i \in \mathbf{I_0},\ j \in \mathbf{I},\ i \neq j \quad (42)$$

$$A_i + fst_i + t_{i0} \leq D_{\max} \quad i \in \mathbf{I} \quad (43)$$

$$A_i \geq 0 \quad i \in \mathbf{I} \quad (44)$$

In the first two constraints, $M$ stands for a very large number. It can be substituted by $D_{\max} + \max_{(i,j) \in \mathbf{I}} \{fst_i + t_{ij}\}$ to avoid out-of-scale values. These constraints state that if a vehicle travels from $i$ to $j$, then $A_j$ should be equal to $(A_i + fst_i + t_{ij})$. The third constraint enforces the *MTD* for each and every customer $i \in \mathbf{I}$. By convention, the arrival time $A_0$ at the depot is fixed to zero.

# 3. Experimental setup and computational results

In the experimentation stage of our study, we implemented an incomplete factorial design rather than a full factorial design of experiments to observe the effect of each independent variable on the solution accuracy and speed. In the sequel, we first present our independent variables and propose an evaluation scheme based on five metrics (key performance indicators). Secondly, we provide the attributes of all 121 test instances making up our test bed, and reveal the hardware and software specs of the employed testing platform. Thirdly, we report our benchmark results in detail. Finally, we propose a distilled guidance as to which combination of mathematical model, valid inequalities, and multithreading options should be used for which size of problems in order to achieve the desired solution accuracy and efficiency.

**Table 1.** List of the factors in the experiments

| No. | Factor | Levels | |
|---|---|---|---|
| 1 | Formulation (mathematical model) | i.<br>ii.<br>iii.<br>iv. | GG<br>MTZ-L<br>BHM<br>MCF |
| 2 | Minimum Number of Vehicles Constraint (MinNV) for GG, MTZ-L and MCF | i.<br>ii. | MinNV(−): constraint omitted<br>MinNV(+): constraint added |
| 3 | Fixed Number of Vehicles Constraint (FixedK) for BHM | i.<br>ii. | FixedK(−): constraint omitted<br>FixedK(+): constraint added |
| 4 | Maximum Number of Vehicles Constraint (MaxNV) for GG, MTZ-L and MCF | i.<br>ii. | MaxNV(−): constraint omitted<br>MaxNV(+): constraint added |
| 5 | Valid Inequalities Configuration (VI-Con) for GG and MTZ-L | i.<br>ii.<br>⋮<br>viii. | VI-000<br>VI-001<br>⋮<br>VI-111 |
| 6 | `Concurrentmip` | i.<br>ii. | 1<br>2 |
| 7 | Valid equalities for MCF | i.<br>ii. | FGX(−) : constraint omitted<br>FGX(+) : constraint added |

## 3.1. Independent variables and KPIs

We basically investigate the effects of seven independent variables (referred to as factors henceforward) as tabulated in Table 1. The fifth factor VI-Con picks one of the eight levels shown in the table which encompass all possible $2^3 = 8$ combinations of the valid inequalities VI-1, VI-2 and VI-3 given in (35), (38) and (39), respectively. In other words, VI-000 implies the configuration free of VI-1, VI-2 and VI-3, while VI-111 points to the configuration with all three VIs included. The sixth factor is the value of the Gurobi solver option `Concurrentmip`, which becomes void unless Gurobi's multithreading/parallel optimization capabilities are turned on. The seventh factor is the inclusion or exclusion of the valid inequalities FGX in (40) in the multi-commodity formulation MCF.

We certainly refrain from a full factorial design of experiments; this would require us to test $(2^3 \times 8) + (2^3 \times 8) + (2^2) + (2^4) = 148$ different experimental conditions for each and every test instance modeled with GG, MTZ-L, BHM or MCF, respectively. Taking into account the size of our test bed, the full factorial approach clearly proves impracticable. Instead, we resort to an incomplete factorial design of experiments, and progress in a sequential approach in which the design of the next scenario of experiments is shaped by the outcome of the current scenario. In order to benchmark our experimental results we propose an evaluation scheme based on the following five key performance indicators (KPIs):

i. CPU time (*t*):

This is the time resource used for solving a given model. Lower *t* is preferred for better efficiency (solution speed). In case the model does not return an optimal solution within a specified time limit (three hours or 10,800 seconds), *t* should be supplemented by other KPIs.

ii. Number of best known solutions found (***#BKS***):

Conforming to the convention in the literature, we indicate with the acronym *BKS* either the best solution reported to date or the objective value of that best solution. The second metric *#BKS* reveals how many best known solutions (*BKS*s) are attained by the tested combination of factors. Each scenario of experiments (tests) is performed on a particular subset of instances. Hence, the target *#BKS* differs from one scenario of experiments to another scenario.

iii. Relative gap between the best feasible and best known objective values (***%BKS***):

Let *BFS* denote the best feasible solution found at the end of a GAMS solver run performed with a given CVRP model. It is actually equal to the tightest (smallest) upper bound on the objective value of that solution satisfying all constraints of the CVRP model including integrality constraints. The third KPI is then computed as *%BKS* = (*BFS*−*BKS*)/*BKS*×100%, which indicates the deviation of *BFS* from the target *BKS*. Lower *%BKS* is preferred for better solution accuracy.

iv. Optimality gap (***%Opt***):

Let the tightest (greatest) lower bound on the objective value of a given CVRP model be denoted by *BPS*, which stands for *best possible solution*. The solution with the objective value equal to *BPS* does not satisfy all integrality constraints unless the solver terminates at proven optimality. Given *BPS* and *BFS*, *%Opt* is computed by GAMS as (|*BFS*−*BPS*|)/|*BFS*|×100% (GAMS Solver Manuals, 2013). It is reported automatically at the end of each run of the solver. In a proven optimal CVRP solution *%Opt* will be zero. Smaller *%Opt* is preferred since it indicates higher accuracy. If a model run completes with zero *%Opt*, one can infer that the solver was able to reach proven optimality within the given time limit.

v. Number of proven optimal solutions found (***#Opt***):

The last metric *#Opt* reveals how many proven optimal solutions are found by the tested combination of factors. The completion of the model run of a given instance with zero *%Opt* increases both metrics *#Opt* and *#BKS* by one.

## 3.2. Test bed and testing platform

### 3.2.1. Characteristics of the test problems

We have set the maximum problem size to 100 customers, and retrieved all CVRPs available from the following online libraries.

i.) **VRP Web** (2015) maintained by NEO (Networking and Emerging Optimization Research Group) at the University of Malaga, Spain.

ii.) **CVRPLIB** (2015) maintained by GALGOS (Algorithms, Optimization and Simulation Group) at the Pontifícia Universidade Católica do Rio de Janeiro, Brasil.

Our test bed accommodates 121 artificially generated symmetric CVRP instances in total. Note that we use the terms "instance" and "test problem" interchangeably. Table A.1 in Appendix A serves as a repository of all tested CVRPs. The column headers of the table show the following attributes for each instance (from left to right): Name, the number of customers ($n$), the uniform vehicle capacity ($Q$), either the number of available vehicles ($K_{\text{fixed}}$) or the number of vehicle routes or tours ($K$) in the *BKS*, tightness which is defined as $\sum_i q_i / (Q \times K_{\text{fixed}})$, i.e. the ratio of the total demand to the total vehicle capacity for a given $K_{\text{fixed}}$, and finally either the proven optimal objective value (*Opt*) or *BKS*. For 118 out of 121 instances, *BKS* is equivalent to *Opt*. This means that the proven optimal objective value is unknown for only three instances. All 121 instances have been sorted in non-decreasing order of $n$ and listed in Table A.1. The interested reader is referred to a recent open-access article by Uchoa et al. (2014) for the origin of each problem.

Sixteen out of 121 instances use two-dimensional Euclidean distances computed with floating-point precision (i.e., without rounding). Moreover, $K$ is unlimited in these 16 instances. In the remaining 104, two-dimensional Euclidean distances are rounded to the nearest integer, and $K$ is fixed. There is only one exception to this, namely, our last test instance X-n101-k25 generated by Uchoa et al. (2014) without fixing $K$. As can be seen in Table A.2, our repository hosts only four time-restricted problems where the traveling time of each route is constrained by a maximum tour duration ($D_{\max}$), and there exists a fixed service time (*fst*) at each customer. Vehicles are assumed to travel at a constant unit speed between any two nodes.

We remark that the *BKS* values of five instances on VRP Web are incorrect or out-of-date. The *Opt* values we obtained from the respective model runs for four of them do not match the *BKS* of VRP Web (see Table A.3). Solution reports and visualizations available from the CVRPLIB website validate this situation as well. One such instance is P-n55-k15. Actually, none of our models was able to return a feasible solution for this instance at the end of three CPU hours. Therefore, we raised the value of $K$ to 16, and substituted this instance with P-n55-k16. Its optimal *BKS* is given in Table A.1. A similar case arose in X-n101-k25. In the proven optimal solution of this instance (Uchoa et al., 2014), there are 25 routes. However, all the feasible solutions we obtained at the end of three CPU hours had 27 routes no matter whether we turned on the MaxNV constraint or not.

We classified the contents of our test bed in non-decreasing order of *n*. This yields the following classes of test instances.

**Class 1**: 39 instances where $n \leq 40$.

**Class 2**: 26 instances where $41 \leq n \leq 50$.

**Class 3**: 29 instances where $51 \leq n \leq 71$.

**Class 4**: 27 instances where $75 \leq n \leq 100$.

**Table 2.** List of GAMS specific and solver specific options applied to all runs.

| GAMS Specific Options and Model Attributes | GUROBI Specific Options (`gurobi.opt` file) |
|---|---|
| `MIP            = GUROBI`<br>`OPTCR          = 0.000`<br>`RESLIM         = 10800`<br>`ITERLIM        = 1.0e9`<br>`ModelName.optfile   = 1`<br>`ModelName.scaleopt  = 1` | `Nodelimit      5.0e8`<br>`Nodefilestart  22.5`<br>`Threads        0`<br>`Concurrentmip  1   or  Concurrentmip  2` |

### 3.2.2. Technical specifications of the testing platform

We performed all our experiments on a Dell Precision T3500 model PC equipped with one Intel Xeon® W3690 3.46 GHz processor and 24 GBytes of random access memory (RAM) of type ECC DDR3. The hexa-core Xeon® processor provides 12 threads with the hyperthreading feature turned on (Intel, 2015). We used the modeling and optimization suite GAMS 24.2.1 to construct and run our CVRP models. The 64-bit Windows version of this software was released on December 9, 2013 by GAMS Corp. Among full-version MILP solvers available under this suite, we chose Gurobi 5.6.0. Table 2 shows the GAMS options, model attributes, and the contents of the file `gurobi.opt` used to set solver-specific parameters (options) in all GAMS runs. `MIP` specifies which solver is to be used for MILP models. By setting `OPTCR` to zero, `RESLIM` to 10,800 and `ITERLIM` to 1.0e9, we run our models until proven optimality, the end of three CPU hours or one billion iterations, whichever occurs first. `ModelName.optfile` tells the solver to read an option file from the GAMS project directory before starting the solution process. `ModelName.scaleopt` turns on the scaling feature of GAMS which applies to all decision variable coefficients and parameter values in the model.

The options in `gurobi.opt` regulate termination conditions and multithreading (parallel optimization) capabilities of Gurobi. `Nodelimit` limits the number of nodes to be explored in the branch-and-bound tree. `Nodefilestart` is the nodefile starting indicator. Recall that our computer has a total memory of 24 GBytes. By setting this option to 22.5 we stop the solver when the node storage in the memory exceeds 22.5 GBytes.

In this way, we reserve 1.5 GBytes memory for the operating system Windows 7, and prevent Gurobi from performing time-costly read and write operations on the hard disk. The options `Threads` and `Concurrentmip` turn on the multithreading (concurrent optimization) capabilities of the MILP solver. When `Threads` is set to zero, the computing load is distributed onto all available six cores (12 threads) of the processor. On the other hand, when `Concurrentmip` is set two, the solver divides available threads evenly between two independent MILP solve operations, and performs them in parallel. Optimization terminates when the first solve operation completes. For an in-depth discussion of the GAMS options, model attributes and Gurobi options, the reader is referred to GAMS Solver Manuals (2013).

We deliberately refrain from tuning the other Gurobi options. The goal in this study is to find a high-performance combination of a mathematical model, valid inequalities, and a generic set of solver options in common that can be applied to all CVRPs as is. A large number of solver options fine-tuned for a particular CVRP instance are not guaranteed to work well in another instance. Moreover, the search space for the best possible problem-dependent configuration of solver options is clearly too big to explore given the sheer number of available options. Therefore we keep all the remaining options at their default values.

### 3.3. Sequential search for the best combination

In this section we develop a sequential funnel approach to be used in our search for the best combination of the factors which were introduced in Table 1. We explain the setup and discuss the results of each scenario of experiments. In our funnel approach, the interpretation of the results obtained from the current scenario sheds light on how the next scenario should be set up. Fig. 1 illustrates a blueprint for this approach. The left-hand side of the depicted funnel tells about which instances have been tested in each scenario. The right-hand side tells about which factors have been fixed and fed into the experiments as input parameters. Finally, the down arrow callout reveals which factor levels are to be adopted before proceeding to the next scenario of experiments. This simple funnel approach helps to narrow down the search for the most effective and efficient combination of factors that is expected to produce the best KPIs. The last scenario in Fig. 1 (Scenario-7) is purposely disconnected from the preceding scenarios. The experiments in that last scenario have been conducted separately to investigate the effect of multithreading / parallel optimization which was not cited as a major factor in Table 1.

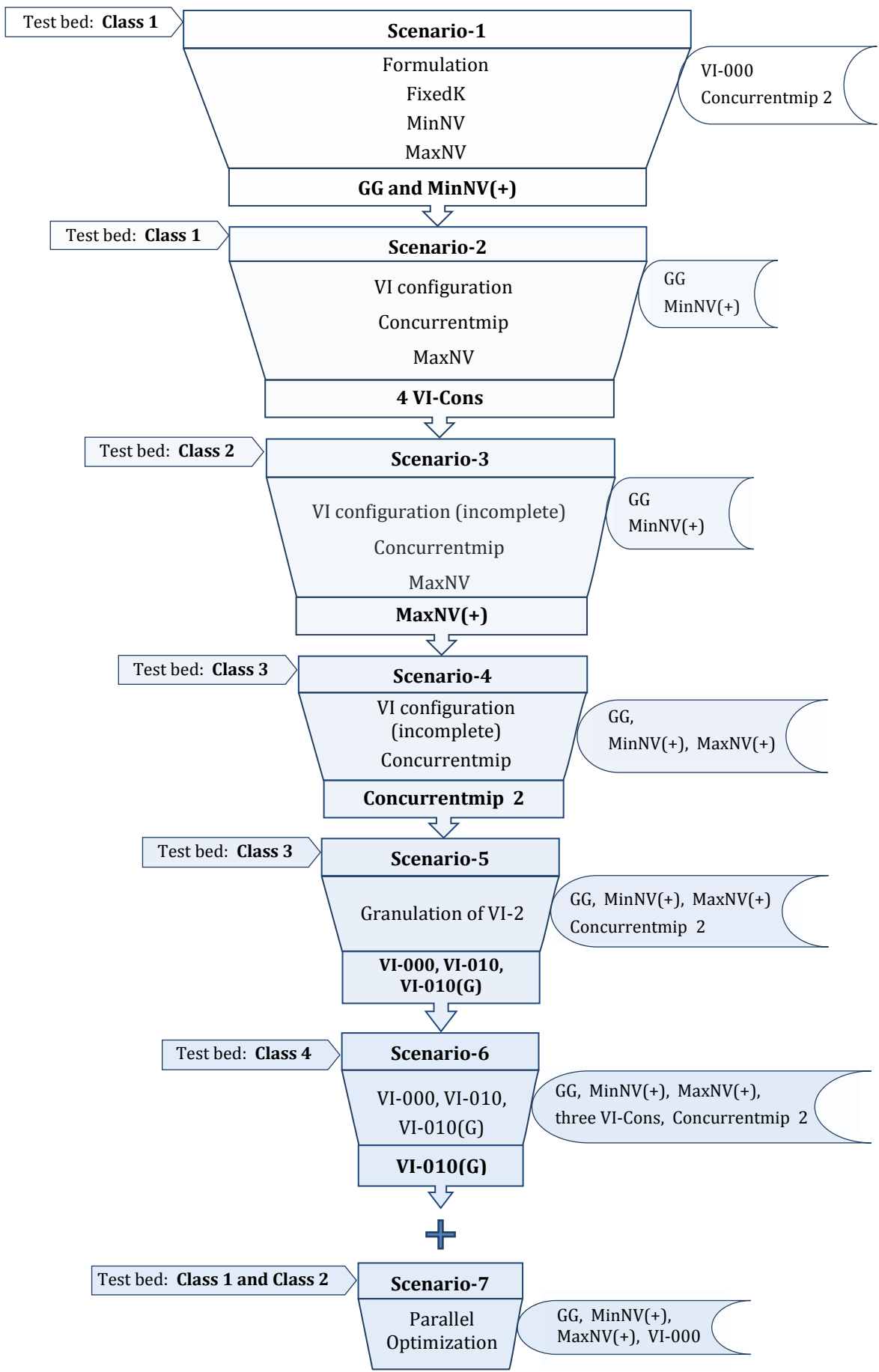

**Fig. 1.** Blueprint for our sequential search

### 3.3.1. Scenario-1: Factors 1, 2, 3 and 4 on Class 1

Our experimentation is initiated with Scenario-1 in which we tested the effect of the first four factors using Class 1 instances. The remaining factors, namely VI-Con and `Concurrentmip`, have been fixed to VI-000 and 2, respectively. The averages of *t*, *%BKS* and *%Opt* along with the *#BKS* and *#Opt* values attained by each combination of factor levels are presented in Table 3.a through Table 3.e. **Boldfaced** figures in these and the subsequent tables indicate **best** KPIs (lowest values of *t*, *%BKS* and *%Opt*; highest values of *#BKS* and *#Opt*). Note that Table 3 shows the test results only for the formulations GG and MTZ-L (the first two levels of factor 1). KPIs for BHM (the third level of factor 1) obtained with FixedK(−) and FixedK(+) and for MCF (the fourth level of factor 1) obtained with MinNV(−) and MinNV(+) are reported and compared to those for GG in Table 4.

**Table 3.a**  Scenario-1 results, part one.

| Average CPU time *t* (s) | MaxNV(−) | MaxNV(+) |
|---|---|---|
| GG | 2053.4 | 724.8 |
| MinNV(−) | 3174.9 | 807.0 |
| MinNV(+) | 931.9 | **642.5** |
| MTZ-L | 7351.4 | 7263.3 |
| MinNV(−) | 8417.8 | 8232.3 |
| MinNV(+) | 6285.0 | 6294.4 |

**Table 3.b**  Scenario-1 results, part two.

| *#BKS* | MaxNV(−) | MaxNV(+) |
|---|---|---|
| GG | 77 | 78 |
| MinNV(−) | 38 | **39** |
| MinNV(+) | **39** | **39** |
| MTZ-L | 49 | 51 |
| MinNV(−) | 23 | 25 |
| MinNV(+) | 26 | 26 |

**Table 3.c**  Scenario-1 results, part three.

| Average *%BKS* | MaxNV(−) | MaxNV(+) |
|---|---|---|
| GG | 0.00% | 0.00% |
| MinNV(−) | **0.00%** | **0.00%** |
| MinNV(+) | **0.00%** | **0.00%** |
| MTZ-L | 0.39% | 0.42% |
| MinNV(−) | 0.48% | 0.44% |
| MinNV(+) | 0.30% | 0.40% |

**Table 3.d**  Scenario-1 results, part four.

| Average *%Opt* | MaxNV(−) | MaxNV(+) |
|---|---|---|
| GG | 0.69% | 0.04% |
| MinNV(−) | 1.33% | 0.07% |
| MinNV(+) | 0.05% | **0.01%** |
| MTZ-L | 17.20% | 17.27% |
| MinNV(−) | 23.79% | 23.90% |
| MinNV(+) | 10.61% | 10.65% |

**Table 3.e**  Scenario-1 results, part five.

| *#Opt* | MaxNV(−) | MaxNV(+) |
|---|---|---|
| GG | 64 | 76 |
| MinNV(−) | 28 | **38** |
| MinNV(+) | 36 | **38** |
| MTZ-L | 30 | 30 |
| MinNV(−) | 12 | 12 |
| MinNV(+) | 18 | 18 |

**Table 4.**   Comparison of the three arc-based formulations BHM, GG and MCF.

| KPI values of BHM, MCF and GG formulations on Class 1 | | | | | |
|---|---|---|---|---|---|
| | BHM Formulation | | MCF Formulation with FGX(−)[a] and MaxNV(−)[b] | | GG Formulation with MinNV(+) and MaxNV(+) |
| | FixedK(−) | FixedK(+) | MinNV(−) | MinNV(+) | |
| *#BKS* | **39** | 36 | 15 | 17 | **39** |
| Avg. *%BKS* | **0.00%** | 0.02% | 22,81% | 20,57% | **0.00%** |
| Max. *%BKS* | **0.00%** | 0.42% | 392,40% | 392,40% | **0.00%** |
| *#Opt* | 32 | 31 | 13 | 14 | **38** |
| Avg. *%Opt* | 0.31% | 0.44% | 15,84% | 14,32% | **0.01%** |
| Max. *%Opt* | 3.02% | 5.06% | 100,00% | 100,00% | **0.45%** |
| Avg. *t* (s) | 2353.4 | 2349.9 | 7611,9 | 7509,4 | **642.5** |
| Median *t* (s) | 38.6 | 59.4 | 10800,3 | 10800,3 | **29.9** |

[a, b]: Computational results of 39 instances of Class 1 show that the KPIs of the MCF formulation improve when there is no coupling between F,G and X and the MaxNV constraint is turned off.

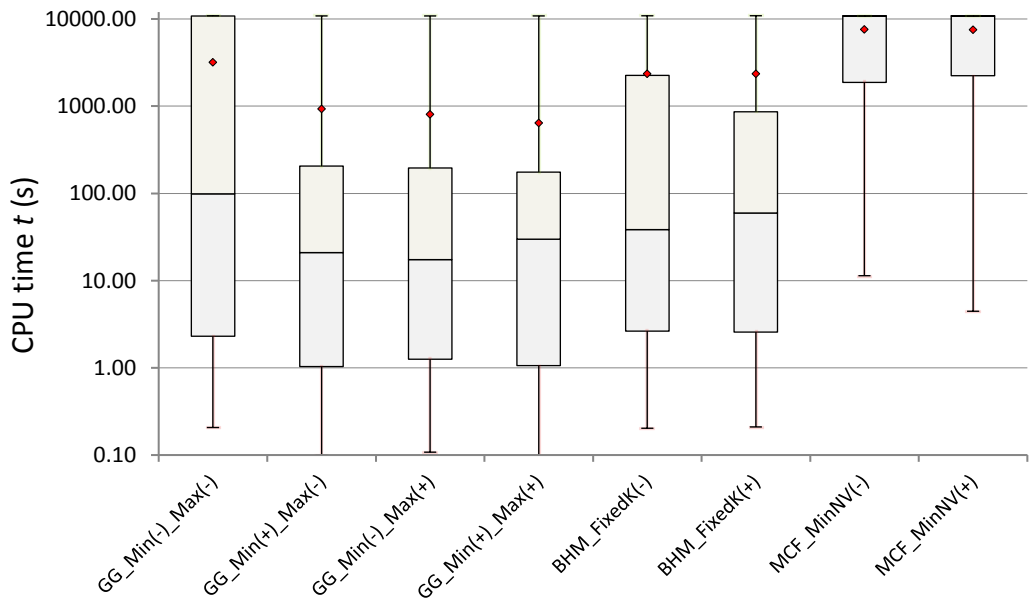

**Fig. 2.**   Box-and-whisker plot of the GG, BHM and MCF solution times for Class 1 instances.

Scenario-1 results strongly suggest that MTZ-L formulation cannot compete with the arc-based formulations GG and BHM even in Class 1 instances. Neither is the MCF formulation a competitor for GG or for BHM. The box-and-whisker plot in Fig. 2 provides a more descriptive look into the solution times of GG, BHM and MCF runs. The ordinate axis of Fig. 2 is displayed in $\log_{10}$ scale. The diamond shaped markers in the chart represent average *t* values. The top and bottom lines of the box are the first and third quartiles, while the band inside the box is the second quartile (i.e., the median) of the plotted data. The upper and lower lines extending vertically from the box are the whiskers. The upper and lower whiskers represent the minimum and maximum points of the data.

We draw two conclusions from Scenario-1: first, GG is superior not only to MTZ-L, but also to BHM and MCF; second, factors 2 and 4 at levels MinNV(+) and MaxNV(+) have a positive effect on the speed of GG. Consequently we adopt GG and MinNV(+) as the levels of choice for factors 1 and 2, and proceed to Scenario-2.

### 3.3.2. Scenario-2: Factors 4, 5 and 6 on Class 1

Scenario-2 is an $8\times2^2$ full factorial design of experiments performed on Class 1 instances. We tested factors 4, 5 and 6 by fixing the first two factors to GG and MinNV(+). CPU time $t$ and *%Opt* are the two KPIs we gauge in this scenario. Their average values attained by each combination of factor levels are presented in Table 5. The best average times in each $t$ column of Table 5 have been shown in *italic*. The overall best average time (384.0) and best *%Opt* (0.00%) are in **bold**. The best grand average $t$ and *%Opt* values in the bottom row have also been marked in **bold**. The grand average values suggest the superiority of MaxNV(+) over MaxNV(–). We observe that CPU time $t$ attains a lower (thus better) average value with MaxNV(+) than with MaxNV(–) in all eight levels of VI-Con and both levels of Concurrentmip. The only exception to this observation occurs in the combination of VI-110 and Concurrentmip 2. A similar situation is observed with *%Opt*. Note that there is a total of 16 distinct combinations of VI-Con and Concurrentmip levels. In twelve of these, the average *%Opt* obtained with MaxNV(+) is smaller than that obtained with MaxNV(–).

**Table 5.** Scenario-2 results.

| Average CPU time $t$ (s), Average *%Opt*, frequency of $t_{min}$ on Class 1 | | | | | | | | | | | | | |
|---|---|---|---|---|---|---|---|---|---|---|---|---|---|
| | MaxNV(–) | | | | | | MaxNV(+) | | | | | | Cum. Freq. of $t_{min}$ |
| VI-Con | Concurrentmip 1 | | | Concurrentmip 2 | | | Concurrentmip 1 | | | Concurrentmip 2 | | | |
| VI-000 | *723.4* | 0.09% | 13 | 931.9 | 0.05% | 9 | 683.4 | 0.08% | 9 | 642.5 | 0.01% | 9 | 40 |
| VI-001 | 1093.7 | 0.08% | 1 | 1154.8 | 0.07% | 0 | 1063.8 | 0.06% | 0 | 857.5 | 0.06% | 2 | 3 |
| VI-010 | 846.0 | *0.03%* | 4 | 893.3 | 0.05% | 7 | **384.0** | **0.00%** | 10 | 676.7 | *0.01%* | 13 | 34 |
| VI-011 | 1159.0 | 0.08% | 0 | 1261.3 | 0.08% | 0 | 1066.5 | 0.09% | 1 | 1123.2 | 0.08% | 1 | 2 |
| VI-100 | 1178.4 | 0.11% | 4 | 866.0 | 0.04% | 7 | 745.5 | 0.07% | 7 | *575.1* | 0.02% | 5 | 23 |
| VI-101 | 1340.6 | 0.14% | 2 | 1135.5 | 0.04% | 5 | 1006.7 | 0.05% | 2 | 1014.8 | 0.10% | 3 | 12 |
| VI-110 | 806.3 | 0.06% | 9 | *673.4* | *0.02%* | 9 | 598.7 | 0.02% | 9 | 810.7 | 0.02% | 6 | 33 |
| VI-111 | 1401.8 | 0.08% | 6 | 1144.7 | 0.08% | 2 | 988.0 | 0.05% | 1 | 977.4 | 0.08% | 0 | 9 |
| Grand Avg. | 1068.6 | 0.09% | — | 1007.6 | **0.05%** | — | **817.1** | **0.05%** | — | 834.7 | **0.05%** | — | |

Another finding in Scenario-2 is that four VI-Con levels (namely, VI-000, VI-010, VI-100 and VI-110) lead the other four levels with respect to $t$. Table 5 presents for each combination of VI-Con and Concurrentmip levels the $t_{min}$ frequency. It is equal to the number of instances in which a given combination proves more efficient than the others. Since Scenario-2 is performed on Class 1 instances, $t_{min}$ frequency can be at most 39. When factor 4 is set to MaxNV(+), the highest $t_{min}$ frequencies are achieved by factor 5 levels VI-000, VI-010, VI-100, and VI-110. Cumulative $t_{min}$ frequencies listed in the last column of Table 5 are also the highest in those levels. Therefore, we adopt the configurations VI-000, VI-010, VI-100, and VI-110 before proceeding to Scenario-3.

### 3.3.3. Scenario-3: Factors 4, 5 and 6 on Class 2

Scenario-3 can be seen as an incomplete repetition of Scenario-2 on Class 2 instances. The difference lies in that only four particular levels of factor 5 (VI-Con) were tested instead of eight. Starting with Scenario-3, we report average *%BKS* as well since they become nonzero as the problem size increases beyond 40 customers. Table 6 shows the average values of three KPIs for each level of factor 5 tested on Class 2 instances. We remark that among 26 Class 2 instances only two have no $K_{fixed}$ (CMT-p01 and CMT-p06). The results in Table 6 can be summarized as follows:

- When Concurrentmip is set to 2, the speed of the Gurobi solver improves at MaxNV(−), but it deteriorates at MaxNV(+). MaxNV(+) seems to spend less CPU times than MaxNV(−). The least average *t* (4142.5 seconds) is attained by the combination of MaxNV(+), Concurrentmip 1 and VI-010.
- Regardless of the level of MaxNV, changing Concurrentmip form 1 to 2 improves average *%Opt* and *%BKS* considerably.
- On average, MaxNV(+) accelerates the solution time of the GG model by 9.4% and 3.3% with Concurrentmip 1 and Concurrentmip 2, respectively.
- Further evidence is needed to make a judgment as to which level of VI-Con is the most effective one.

By the end of this scenario, we adopt MaxNV(+) level for factor 4, and proceed to Scenario-4 to further investigate the effect of the last two factors on Class 3 instances.

**Table 6.** Scenario-3 results.

| Average CPU time *t* (s), Average *%BKS*, Average *%Opt* on Class 2 | | | | | | | | | | | | |
|---|---|---|---|---|---|---|---|---|---|---|---|---|
| VI-Con | MaxNV(−) | | | | | | MaxNV(+) | | | | | |
| | Concurrentmip 1 | | | Concurrentmip 2 | | | Concurrentmip 1 | | | Concurrentmip 2 | | |
| VI-000 | 4956.5 | 0.52% | 1.50% | 4923.8 | 0.25% | 1.24% | 4561.0 | 0.59% | 1.49% | 4560.9 | 0.25% | 1.24% |
| VI-010 | 4599.5 | 0.46% | 1.44% | 4744.7 | 0.25% | 1.24% | **4142.5** | 0.41% | 1.39% | 4638.4 | **0.16%** | **1.20%** |
| VI-100 | 4855.5 | 0.74% | 1.70% | 4515.3 | 0.45% | 1.39% | 4757.3 | 0.82% | 1.88% | 4678.5 | 0.49% | 1.50% |
| VI-110 | 5640.1 | 0.59% | 1.57% | 5162.3 | 0.31% | 1.44% | 4703.0 | 0.58% | 1.48% | 4835.6 | 0.43% | 1.44% |
| Grand Avg. | 5012.9 | 0.58% | 1.55% | 4836.5 | **0.32%** | **1.33%** | **4541.0** | 0.60% | 1.56% | 4678.3 | 0.33% | 1.35% |

### 3.3.4. Scenario-4: Factors 5 and 6 on Class 3

Having fixed factor 4 to level MaxNV(+), we move on to Class 3 instances in Scenario-4. In this class the problem size *n* varies between 51 and 71. Levels VI-000, VI-010, VI-100 and VI-110 of factor 5 have been tested on 29 instances for both levels of factor 6. We observe that the ratio of *#BKS* to *n* drops drastically. This makes *#BKS* a pivotal metric worthy of reporting alongside *%BKS* and *%Opt* in Table 7.a. In each tested combination of factors 5 and 6, only four instances reach zero *%Opt* within three hours. Table 7.b reveals the solution times of those four instances only.

**Table 7.a**  Scenario-4 results, part one.

| Average *%BKS*, Average *%Opt*, *#BKS* on Class 3 | | | | | | |
|---|---|---|---|---|---|---|
| VI-Con | Concurrentmip 1 | | | Concurrentmip 2 | | |
| VI-000 | 0.79% | 3.22% | 8 | **0.49%** | **2.92%** | 10 |
| VI-010 | 0.79% | 3.42% | **11** | 0.61% | 3.08% | 9 |
| VI-100 | 0.88% | 3.34% | 10 | 0.75% | 3.18% | 6 |
| VI-110 | 0.76% | 3.36% | 10 | 0.66% | 3.18% | 7 |
| Grand Avg. | 0.80% | 3.33% | **9.75** | **0.63%** | **3.09%** | 8.00 |

**Table 7.b**  Scenario-4 results, part two.

| CPU time *t* (s) for Class 3 solved with zero *%Opt* | | | | | | | | |
|---|---|---|---|---|---|---|---|---|
| | Concurrentmip 1 | | | | Concurrentmip 2 | | | |
| VI-Con | B-n52-k7 | B-n56-k7 | B-n64-k9 | F-n72-k4 | B-n52-k7 | B-n56-k7 | B-n64-k9 | F-n72-k4 |
| VI-000 | **342.9** | 506.5 | 3672.3 | 166.7 | 494.9 | 658.3 | 5076.5 | 142.6 |
| VI-010 | 510.0 | 600.7 | **1523.3** | 174.5 | 524.9 | 680.4 | 8408.7 | 158.7 |
| VI-100 | 399.9 | **389.4** | 2242.9 | **81.7** | 629.6 | 599.1 | 1755.5 | 127.5 |
| VI-110 | 499.9 | 718.1 | 8133.0 | 143.4 | 938.5 | 728.6 | 2920.9 | 107.0 |
| Average *t* | **438.2** | **553.7** | **3892.9** | 141.6 | 647.0 | 666.6 | 4540.4 | **133.9** |

In order to make a more informed decision about factor 6, we utilize the box-and-whisker plots of the *%BKS* and *%Opt* results as shown in Fig. 3. The box heights for `Concurrentmip 2` are smaller than those for `Concurrentmip 1` especially in the plot of *%BKS*. It can be interpreted as a lesser degree of variability in the *%BKS* results, which leads us to adopt `Concurrentmip 2`, and proceed to Scenario-5.

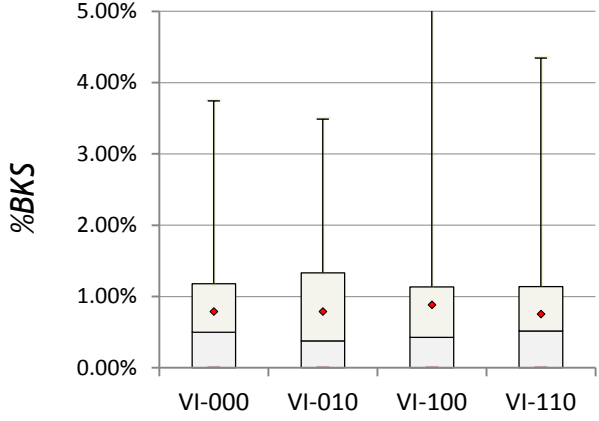

**Fig. 3.a** Class 3 *%BKS* results with `Concurrentmip 1`.

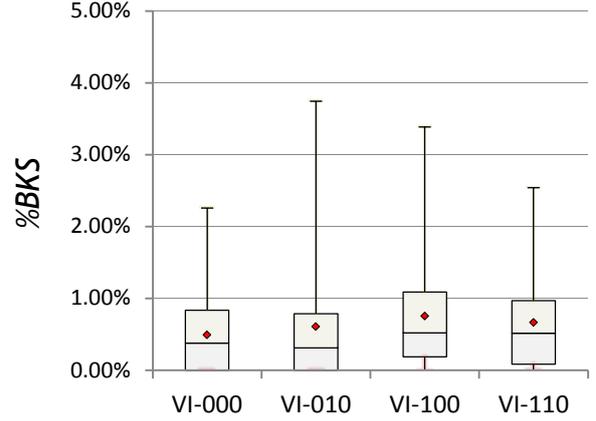

**Fig. 3.b** Class 3 *%BKS* results with `Concurrentmip 2`.

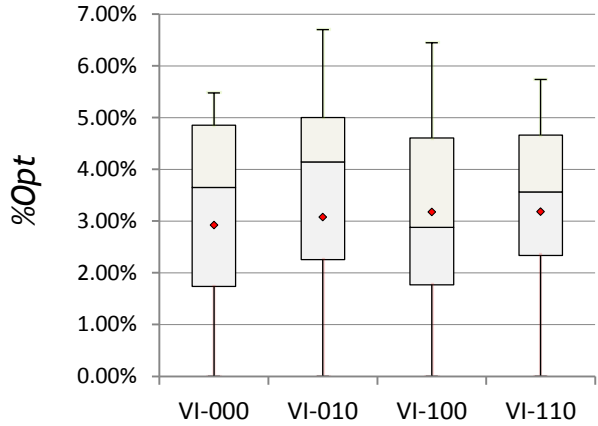

**Fig. 3.c** Class 3 *%Opt* results with `Concurrentmip 1`.

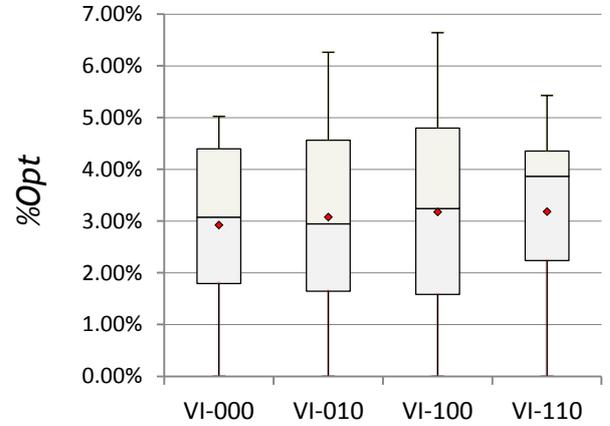

**Fig. 3.d** Class 3 *%Opt* results with `Concurrentmip 2`.

### 3.3.5. Scenario-5: Granulation of VI-2 on Class 3

In this scenario we experiment again with Class 3 instances where we try the so-called granulation of VI-2. Recall that VI-2 in (38) is added to the GG model as an extra set of constraints preventing subtours of size two. The GG model is certainly self-sufficient in preventing subtours of any size. As discussed in Section 2, VI-2 and VI-3 are supplementary constraints serving the same purpose. We design a simple topology-driven granulation scheme for VI-2. Note that the number of VI-2 equations in VI-010 and VI-110 configurations before granulation is equal to $\binom{n}{2}$. We seek to reduce this number, thereby consume less memory space during the solution process, and eventually improve the KPIs. Let the problem space diameter $D$ be defined as $D = \max_{(i,j) \in A} \{d_{ij}\}$. Define a threshold distance parameter $T = \dfrac{D}{\lceil \ln(n+1) \rceil}$. We shall write VI-2 for customer pair $(i, j)$ if and only if $d_{ij} \leq T$ holds true. The underlying idea is if customers $i$ and $j$ are apart from one another by more than $T$, then the formation of the subtour ($i{\rightarrow}j{\rightarrow}i$) is rather unlikely to happen. So, there is no need to supplement the GG model with the VI-2 equation involving these two customers. The GG model itself prevents subtours anyway.

The granular version of VI-2 imposes clearly a lesser number of constraints on the GG model. Table 8 presents the VI-2 counts before and after granulation in Class 3 and Class 4 instances. The proposed granulation scheme achieves an average reduction rate of 80.5% in the number of VI-2 equations.

**Table 8.** Reduction in the count of VI-2 equations through simple granulation.

| Instances | Prob. size $n$ | Average Values | | |
|---|---|---|---|---|
| | | Full Count of VI-2 eqs. | Count of VI-2 eqs. in VI-010(G) | Reduction by |
| Class 3 | 60 | 1779 | 318 | 82.1% |
| Class 4 | 87 | 3842 | 812 | 78.9% |
| Grand Avg. | 73 | 2774 | 556 | 80.5% |

**Table 9.** Scenario-5 results.

| Average KPI values on Class 3 with `Concurrentmip 2` and MaxNV(+) | | | | | |
|---|---|---|---|---|---|
| VI-Con | *%BKS* | *%Opt* | #BKS | #Opt | $t$ (s) for four instances |
| VI-000 | **0.49%** | **2.92%** | 10 | 4 | 1593.1 |
| VI-010 | 0.61% | 3.08% | 9 | 4 | 2443.2 |
| VI-010(G) | **0.49%** | 2.98% | **13** | 4 | 1096.8 |
| VI-100 | 0.75% | 3.18% | 6 | 4 | **777.9** |
| VI-110 | 0.66% | 3.18% | 7 | 4 | 1173.7 |
| VI-110(G) | 0.67% | 3.22% | 9 | 3 | 3072.8 |

The granular versions of factor 5 levels are denoted by VI-010(G) and VI-110(G). Table 9 contains average KPIs attained by the former and newly tested levels of factor 5. The last column shows the average CPU time spent at each level of factor 5 for finding a proven optimal solution for the instances `B-n52-k7`, `B-n56-k7`, `B-n64-k9` and `F-n72-k4`. Recall that only these four instances can be solved in less than three hours. The top three *%BKS* and *%Opt* performers in Table 9 are VI-000, VI-010(G) and VI-010 whereas VI-100, VI-110 and once more VI-010(G) are the three most efficient VI configurations. At this point we prioritized accuracy over speed, and adopted VI-000, VI-010 and VI-010(G) before moving on to Scenario-6.

### 3.3.6. Scenario-6: Three VI configurations on Class 4

In this scenario of experiments of our sequential funnel approach, factors 1, 2, 4 and 6 are fixed to levels GG, MinNV(+), MaxNV(+) and `Concurrentmip 2`, respectively. In factor 5 we eliminated all levels but VI-000 and VI-010, and added VI-010(G) as a new level to consider. In short, we test only three VI configurations on Class 4 instances. Recall from Table A.1 that in Class 4 there are 27 instances in total; 12 with 75, one with 77, one with 79 and 13 with 100 customers. Three of these instances have a maximum tour duration constraint. The results of Scenario-6 are shown in Table 10.a and Fig. 4. At most four instances can be solved in less than three hours with zero optimality gap. A comparison of CPU times of these four is given in Table 10.b. Descriptive statistics gathered from Scenario-6 tell us the following:

- The newly added granular factor 5 level VI-010(G) outperforms VI-000 and VI-010 in average solution accuracy. Although its worst-case *%BKS* (the fourth quartile) is higher than that of VI-010 (5.88% vs. 3.30%), it has a lower median *%BKS* than the other two VI configurations (0.35% vs. 0.74% and 1.04%). We observe a better range of gaps with VI-010(G) also in *%Opt*.

- VI-010(G) attains the best average *t* on the four instances shown in Table 10.b. Among the three VI configurations it is the only one which can solve all four instances to proven optimality.

In the light of the above analysis, we suggest VI-010(G) as the VI configuration of choice for Class 4 and larger instances. We proceed to the last scenario of our experiments.

**Table 10.a**  Scenario-6 results, part one.

| Average KPI values on Class 4 with `Concurrentmip 2` and MaxNV(+) | | | | |
|---|---|---|---|---|
| VI-Con | *%BKS* | *%Opt* | *#BKS* | *#Opt* |
| VI-000 | 1.16% | 4.96% | **8** | 3 |
| VI-010 | 0.96% | 4.70% | 6 | 2 |
| VI-010(G) | **0.86%** | **4.67%** | 6 | **4** |

**Table 10.b**  Scenario-6 results, part two.

| CPU time *t* (s) for Class 4 solved with zero *%Opt* | | | | | |
|---|---|---|---|---|---|
| VI-Con | P-n76-k4 | tail75d | P-n101-k4 | M-n101-k10 | Average |
| VI-000 | 5042.3 | 10806.1 | 2460.1 | **4855.5** | 5791.0 |
| VI-010 | 5041.7 | 10808.7 | **1044.8** | 10806.4 | 6925.4 |
| VI-010(G) | **2160.0** | **7039.1** | 1201.4 | 5608.0 | **4002.1** |

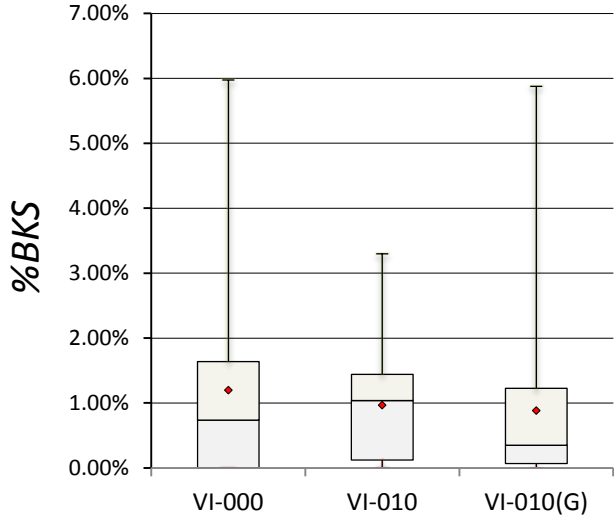
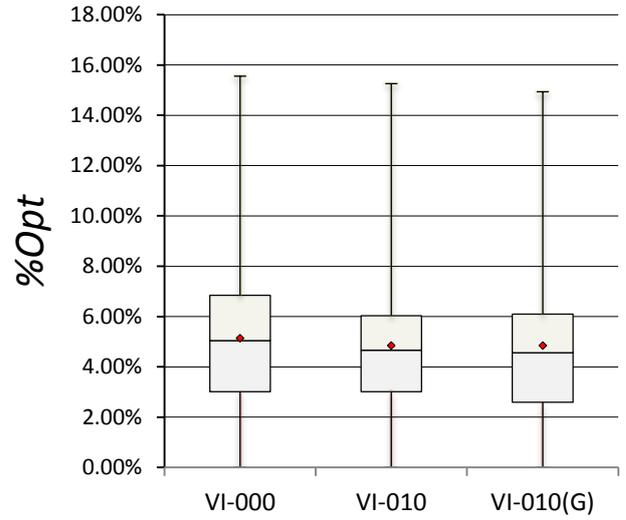

**Fig. 4.a** Class 4 *%BKS* results with `Concurrentmip 2`   **Fig. 4.b** Class 4 *%Opt* results with `Concurrentmip 2`.

### 3.3.7. Scenario-7: Benefit of parallel optimization on Class 1 and Class 2

Scenario-7 has been depicted as disconnected from the funnel in Fig. 1 since it was designed to assess the effect of parallel optimization which was not included in Table 1 as a major factor. We fix factors 1, 2, 4 and 5 to levels GG, MinNV(+), MaxNV(+) and VI-000, respectively. Then we test the following three settings using a test bed comprised of Class 1 and Class 2 which have a combined total of 65 instances:

i. Parallel optimization (multithreading) with `Concurrentmip 1`.
ii. Parallel optimization with `Concurrentmip 2`.
iii. Single-threaded optimization with the solver options `Threads` and `Concurrentmip` turned off.

Table 11.a and Table 11.b present the results of Scenario-7 in detail. The advantage of parallel optimization becomes most apparent in the KPI *t*. The median *t* of single-threaded model runs is around 4.2 and 3.7 times that of multithreaded runs with `Concurrentmip 1` and `Concurrentmip 2`, respectively. In one particular instance, namely P-n50-k8, the maximum *%BKS* and *%Opt* values produced by `Concurrentmip 1` are noticeably worse than those of the single-threaded runs. Yet, it is an exceptional case.

Table 11.b offers a more in-depth comparison between multi- and single-threaded runs in two ways, namely accuracy and speed. First, we separate 18 instances whose single-threaded runs consumed the 3-hour time limit, thereby terminated at nonoptimality (*%Opt* > 0). We report the average *%BKS* and *%Opt* values of the multi- and single-threaded runs of this subset of instances in Table 11.b. The descriptive statistics in the table certify the accuracy superiority of multithreading at the `Concurrentmip 2` level of factor 6. Secondly, we consider the remaining 47 instances of the joint test bed. The single-threaded runs of those instances have been completed with zero *%Opt*, that is to say their solution times have not been cut off at the 3-hour time limit. Hence, the speed of multithreading can be benchmarked equitably on them. Table 11.b reveals the average, median and maximum *t* values measured in three different settings. Turning on multithreading in Gurobi by setting the solver option `Concurrentmip` to 1 (2) reduces the average *t* by 68% (66%) and the median *t* by 64% (73%). We conclude that multithreading is recommended not only for higher accuracy (effectiveness), but also for higher speed (efficiency) in CVRP solutions.

**Table 11.a**  Scenario-7 results, part one.

| KPI values on Class 1 and Class 2 combined with GG, MinNV(+), MaxNV(+) and VI-000 | | | |
|---|---|---|---|
| | Multithreading | | Single Threading |
| | `Concurrentmip 1` | `Concurrentmip 2` | `n/a` |
| *#BKS* | 61 | 61 | 58 |
| Avg. *%BKS* | 0.16% | 0.07% | 0.15% |
| Median *%BKS* | 0.00% | 0.00% | 0.00% |
| Max. *%BKS* | 9.06% | 3.18% | 3.26% |
| *#Opt* | 57 | 55 | 47 |
| Avg. *%Opt* | 0.56% | 0.49% | 0.75% |
| Median *%Opt* | 0.00% | 0.00% | 0.00% |
| Max. *%Opt* | 12.54% | 7.55% | 7.50% |
| Avg. *t* (s) | 1887.4 | 2261.3 | 3649.3 |
| Median *t* (s) | 177.3 | 200.6 | 739.8 |

**Table 11.b**  Scenario-7 results, part two.

| KPI values attained by single-threaded runs of Class 1 and Class 2 | | | |
|---|---|---|---|
| | Multithreading | | Single Threading |
| | `Concurrentmip 1` | `Concurrentmip 2` | `n/a` |
| # of nonoptimal runs | 8 | 10 | 18 |
| Avg. *%BKS* of the instances with nonoptimal single-threaded runs | 0.59% | 0.24% | 0.54% |
| Avg. *%Opt* of the instances with nonoptimal single-threaded runs | 2.01% | 1.76% | 2.70% |
| Avg. *t* (s) of the instances with optimal single-threaded runs | 288.7 | 309.3 | 910.6 |
| Median *t* (s) of the instances with optimal single-threaded runs | 44.9 | 33.3 | 124.9 |
| Max. *t* (s) of the instances with optimal single-threaded runs | 4,601.6 | 2,927.6 | 5,865.6 |

## 3.4. Summary of the experimental results

In Table 12 we present a concise summary of our experimental results reported thus far. With this table we intend to offer a distilled guidance to practitioners working on routing problems and to the CVRP research community at large. The table suggests using the single-commodity arc-based formulation, i.e. the GG model. Constraints bounding the minimum and the maximum number of vehicles dispatched from the depot should be incorporated into the GG model. Multithreading and scaling features of the Gurobi solver should be turned on. We recommend VI-010 and `Concurrentmip 1` for problems with no more than 40 customers. The granular version VI-010(G) and `Concurrentmip 2` are recommended for problems with 75 or more customers. For problems of sizes in between, our recipe branches into two depending on whether accuracy or speed is the number one priority. If the problem at hand is to be solved to proven optimality, i.e., if the target optimality gap is zero, then we recommend to couple VI-010 with `Concurrentmip 1` for shortest CPU times. However, if the problem is most likely unsolvable to optimality within the specified time limit, then we recommend to couple either VI-010 or VI-010(G) with `Concurrentmip 2` for the smallest gap.

The last three columns of Table 12 report the average levels of solution speed and accuracy/bounding performance that can be acquired with the suggested model, solver, VI-Con and `Concurrentmip` combination. The average CPU time value $t$ is deliberately not provided for those instances which have more than 50 customers (Class 3 and Class 4), since only 8 out of 56 such instances can be solved to proven optimality in less than three hours. For the remaining 48 instances the solution process was interrupted after three hours; therefore, their CPU times cannot be compared directly.

**Table 12.** A distilled guidance as the outcome of this study.

| Model | Extra constraints and features | Range of $n$[1] | VI-Con | | Gurobi solver option `Concurrentmip` | Average Performance Levels | | |
|---|---|---|---|---|---|---|---|---|
| | | | | | | CPU time $t$ (s) | %BKS | %Opt |
| GG | ▪ MinNV ▪ MaxNV ▪ Multithreading and scaling | $n \leq 40$ | | VI-010 | 1 | 384.0 | 0.00% | 0.00% |
| | | $41 \leq n \leq 50$ | CPU[2]: | VI-010 | 1 | 4142.5 | 0.41% | 1.39% |
| | | | %gap[3]: | VI-010 | 2 | 4638.4 | 0.16% | 1.20% |
| | | $51 \leq n < 75$ | CPU: | VI-010 | 1 | – | 0.79% | 3.42% |
| | | | %gap: | VI-010(G) | 2 | – | 0.49% | 2.98% |
| | | $n \geq 75$ | | VI-010(G) | 2 | – | 0.86% | 4.67% |

[1] Number of customers in the given CVRP instance.

[2] If CPU time has the priority.

[3] If *%BKS* or *%Opt* has the priority.

## 4. Conclusion and future research

In our experimental study we have gathered a considerable amount of empirical evidence to verify the superiority of the single-commodity and two-commodity arc-based formulations (GG and BHM, respectively) over the multi-commodity arc-based formulation and the node-based formulation (MCF and MTZ-L, respectively). We have observed that the appropriate choice of the following factors can make a substantial difference in both solution accuracy and speed:

i. Mathematical model.
ii. Valid inequalities.
iii. Multithreading (parallel optimization) options of the employed MILP solver.

Our findings suggest that solving the GG model turns out to be a viable approach especially for instances with more than 50 customers. Supplementing the GG model with the constraints that impose the minimum required number of vehicles (MinNV in (32)) and the maximum number of available vehicles where applicable (MaxNV in (33)) improves the average solution speed and reduces the optimality gaps as well. Among all valid inequality configurations taken into account during experimentation, VI-010 which involves only the generalized subtour elimination constraint of size two (constraints (38)) stands out as the most efficient configuration. For test instances with more than 50 customers we tried to granulate VI-2 based on the dispersion of customers in the problem space. Our simple topology-driven granulation technique reduces the number of VI-2 equations to be added to the GG model approximately by 80% across 56 instances. Granulation of VI-2 achieves a favorable trade-off between solution speed and accuracy. Furthermore, turning on the multithreading feature of the Gurobi solver yields a significant improvement in the solution speed of the GG model. Last but not least, we believe that a more detailed analytical analysis of CVRP formulations is still worthwhile. New tailored valid inequalities can be investigated especially for the BHM model which we found to be the second most competitive model after GG.

# Appendix A

**Table A.1.** Attributes of all 121 instances in the test bed sorted in non-decreasing order of $n$.

| Instance Name | $n$ | $Q$ | $K_{fixed}$ / $K$ in BKS | Tightness | Opt or BKS |
|---|---|---|---|---|---|
| **Class 1**: 39 instances with $n \leq 40$ | | | | | |
| E-n13-k4 | 12 | 6,000 | 4 / – | 0.76 | 247* |
| P-n16-k8 | 15 | 35 | 8 / – | 0.88 | 450* |
| ulysses-n16-k3 | 15 | 5 | 3 / – | 1.00 | 8,232* |
| gr-n17-k3 | 16 | 6 | 3 / – | 0.89 | 2,685* |
| P-n19-k2 | 18 | 160 | 2 / – | 0.97 | 212* |
| P-n20-k2 | 19 | 160 | 2 / – | 0.97 | 216* |
| gr-n21-k3 | 20 | 7 | 3 / – | 0.95 | 3,704* |
| P-n21-k2 | 20 | 160 | 2 / – | 0.93 | 211* |
| E-n22-k4 | 21 | 6,000 | 4 / – | 0.94 | 375* |
| P-n22-k2 | 21 | 160 | 2 / – | 0.96 | 216* |
| P-n22-k8 | 21 | 3,000 | 8 / – | 0.94 | 603* |
| ulysses-n22-k4 | 21 | 6 | 4 / – | 0.88 | 9,312* |
| E-n23-k3 | 22 | 4,500 | 3 / – | 0.75 | 569* |
| P-n23-k8 | 22 | 40 | 8 / – | 0.98 | 529* |
| gr-n24-k4 | 23 | 7 | 4 / – | 0.82 | 2,053* |
| fri-n26-k3 | 25 | 10 | 3 / – | 0.83 | 1,353* |
| bayg-n29-k4 | 28 | 8 | 4 / – | 0.88 | 2,050* |
| bays-n29-k5 | 28 | 6 | 5 / – | 0.93 | 2,963* |
| E-n30-k3 | 29 | 4,500 | 3 / – | 0.94 | 534* |
| E-n30-k4 | 29 | 4,500 | 4 / – | 0.71 | 503* |
| B-n31-k5 | 30 | 100 | 5 / – | 0.82 | 672* |
| E-n31-k7 | 30 | 140 | 7 / – | 0.92 | 379* |
| A-n32-k5 | 31 | 100 | 5 / – | 0.82 | 784* |
| A-n33-k5 | 32 | 100 | 5 / – | 0.89 | 661* |
| A-n33-k6 | 32 | 100 | 6 / – | 0.90 | 742* |
| E-n33-k4 | 32 | 8,000 | 4 / – | 0.92 | 835* |

| Instance Name | n | Q | K_fixed / K in BKS | Tightness | Opt or BKS |
|---|---|---|---|---|---|
| A-n34-k5 | 33 | 100 | 5 / – | 0.92 | 778* |
| B-n34-k5 | 33 | 100 | 5 / – | 0.91 | 788* |
| B-n35-k5 | 34 | 100 | 5 / – | 0.87 | 955* |
| A-n36-k5 | 35 | 100 | 5 / – | 0.88 | 799* |
| A-n37-k5 | 36 | 100 | 5 / – | 0.81 | 669* |
| A-n37-k6 | 36 | 100 | 6 / – | 0.95 | 949* |
| A-n38-k5 | 37 | 100 | 5 / – | 0.96 | 730* |
| B-n38-k6 | 37 | 100 | 6 / – | 0.85 | 805* |
| A-n39-k5 | 38 | 100 | 5 / – | 0.95 | 822* |
| A-n39-k6 | 38 | 100 | 6 / – | 0.88 | 831* |
| B-n39-k5 | 38 | 100 | 5 / – | 0.88 | 549* |
| P-n40-k5 | 39 | 140 | 5 / – | 0.88 | 458* |
| B-n41-k6 | 40 | 100 | 6 / – | 0.95 | 829* |
| **Class 2**: 26 instances with $41 \leq n \leq 50$ | | | | | |
| dantzig-n42-k4 | 41 | 11 | 4 / – | 0.93 | 1,142* |
| swiss-n42-k5 | 41 | 9 | 5 / – | 0.91 | 1,668* |
| B-n43-k6 | 42 | 100 | 6 / – | 0.87 | 742* |
| A-n44-k6 | 43 | 100 | 6 / – | 0.95 | 937* |
| B-n44-k7 | 43 | 100 | 7 / – | 0.92 | 909* |
| A-n45-k6 | 44 | 100 | 6 / – | 0.99 | 944* |
| A-n45-k7 | 44 | 100 | 7 / – | 0.91 | 1,146* |
| B-n45-k5 | 44 | 100 | 5 / – | 0.97 | 751* |
| B-n45-k6 | 44 | 100 | 6 / – | 0.99 | 678* |
| F-n45-k4 | 44 | 2,010 | 4 / – | 0.90 | 724* |
| P-n45-k5 | 44 | 150 | 5 / – | 0.92 | 510* |
| A-n46-k7 | 45 | 100 | 7 / – | 0.86 | 914* |
| A-n48-k7 | 47 | 100 | 7 / – | 0.89 | 1,073* |
| att-n48-k4 | 47 | 15 | 4 / – | 0.78 | 40,002* |
| gr-n48-k3 | 47 | 16 | 3 / – | 0.98 | 5,985* |
| hk-n48-k4 | 47 | 15 | 4 / – | 0.78 | 14,749* |
| B-n50-k7 | 49 | 100 | 7 / – | 0.87 | 741* |
| B-n50-k8 | 49 | 100 | 8 / – | 0.92 | 1,312* |
| P-n50-k10 | 49 | 100 | 10 / – | 0.95 | 696* |
| P-n50-k7 | 49 | 150 | 7 / – | 0.91 | 554* |
| P-n50-k8 | 49 | 120 | 8 / – | 0.99 | 631* |
| B-n51-k7 | 50 | 100 | 7 / – | 0.98 | 1,032* |
| CMT-p01 | 50 | 160 | – / 5 | n/a | 524.61* |
| CMT-p06 | 50 | 160 | – / 6 | n/a | 555.43 |
| E-n51-k5 | 50 | 160 | 5 / – | 0.97 | 521* |
| P-n51-k10 | 50 | 80 | 10 / – | 0.97 | 741* |
| **Class 3**: 29 instances with $51 \leq n \leq 71$ | | | | | |
| B-n52-k7 | 51 | 100 | 7 / – | 0.87 | 747* |
| A-n53-k7 | 52 | 100 | 7 / – | 0.95 | 1,010* |
| A-n54-k7 | 53 | 100 | 7 / – | 0.96 | 1,167* |
| A-n55-k9 | 54 | 100 | 9 / – | 0.93 | 1,073* |
| P-n55-k10 | 54 | 115 | 10 / – | 0.91 | 694* |
| P-n55-k15 (substituted with P-n55-k16) | 54 | 70 | 15 / – | 0.99 | 989* |
| P-n55-k7 | 54 | 170 | 7 / – | 0.88 | 568* |
| P-n55-k8 | 54 | 160 | 8 / – | 0.81 | 588* |
| B-n56-k7 | 55 | 100 | 7 / – | 0.88 | 707* |
| B-n57-k7 | 56 | 100 | 7 / – | 1.00 | 1,153* |
| B-n57-k9 | 56 | 100 | 9 / – | 0.89 | 1,598* |
| A-n60-k9 | 59 | 100 | 9 / – | 0.92 | 1,354* |
| P-n60-k10 | 59 | 120 | 10 / – | 0.95 | 744* |

| Instance Name | n | Q | K_fixed / K in BKS | Tightness | Opt or BKS |
|---|---|---|---|---|---|
| P-n60-k15 | 59 | 80 | 15 / – | 0.95 | 968* |
| A-n61-k9 | 60 | 100 | 9 / – | 0.98 | 1,034* |
| A-n62-k8 | 61 | 100 | 8 / – | 0.92 | 1,288* |
| A-n63-k10 | 62 | 100 | 10 / – | 0.93 | 1,314* |
| A-n63-k9 | 62 | 100 | 9 / – | 0.97 | 1,616* |
| B-n63-k10 | 62 | 100 | 10 / – | 0.92 | 1,496* |
| A-n64-k9 | 63 | 100 | 9 / – | 0.94 | 1,401* |
| B-n64-k9 | 63 | 100 | 9 / – | 0.98 | 861* |
| A-n65-k9 | 64 | 100 | 9 / – | 0.97 | 1,174* |
| P-n65-k10 | 64 | 130 | 10 / – | 0.94 | 792* |
| B-n66-k9 | 65 | 100 | 9 / – | 0.96 | 1,316* |
| B-n67-k10 | 66 | 100 | 10 / – | 0.91 | 1,032* |
| B-n68-k9 | 67 | 100 | 9 / – | 0.93 | 1,272* |
| A-n69-k9 | 68 | 100 | 9 / – | 0.94 | 1,159* |
| P-n70-k10 | 69 | 135 | 10 / – | 0.97 | 827* |
| F-n72-k4 | 71 | 30,000 | 4 / – | 0.96 | 237* |
| **Class 4**: 27 instances with $75 \leq n \leq 100$ | | | | | |
| CMT-p02 | 75 | 140 | – / 10 | n/a | 835.26* |
| CMT-p07 | 75 | 140 | – / 11 | n/a | 909.68 |
| E-n76-k10 | 75 | 140 | 10 / – | 0.97 | 830* |
| E-n76-k14 | 75 | 100 | 14 / – | 0.97 | 1,021* |
| E-n76-k7 | 75 | 220 | 7 / – | 0.89 | 682* |
| E-n76-k8 | 75 | 180 | 8 / – | 0.95 | 735* |
| P-n76-k4 | 75 | 350 | 4 / – | 0.97 | 593* |
| P-n76-k5 | 75 | 280 | 5 / – | 0.97 | 627* |
| tai75a | 75 | 1,445 | – / 10 | n/a | 1,618.36* |
| tai75b | 75 | 1,679 | – / 9 | n/a | 1,344.64* |
| tai75c | 75 | 1,122 | – / 9 | n/a | 1,291.01* |
| tai75d | 75 | 1,699 | – / 9 | n/a | 1,365.42* |
| B-n78-k10 | 77 | 100 | 10 / – | 0.94 | 1,221* |
| A-n80-k10 | 79 | 100 | 10 / – | 0.94 | 1,763* |
| CMT-p03 | 100 | 200 | – / 8 | n/a | 826.14* |
| CMT-p08 | 100 | 200 | – / 9 | n/a | 865.94 |
| CMT-p12 | 100 | 200 | – / 10 | n/a | 819.56* |
| CMT-p14 | 100 | 200 | – / 11 | n/a | 866.37 |
| E-n101-k8 | 100 | 200 | 8 / – | 0.91 | 815* |
| E-n101-k14 | 100 | 112 | 14 / – | 0.93 | 1,067* |
| M-n101-k10 | 100 | 200 | 10 / – | 0.91 | 820* |
| P-n101-k4 | 100 | 400 | 4 / – | 0.91 | 681* |
| tai100a | 100 | 1,409 | – / 11 | n/a | 2,041.34* |
| tai100b | 100 | 1,842 | – / 11 | n/a | 1,939.90* |
| tai100c | 100 | 2,043 | – / 11 | n/a | 1,406.20* |
| tai100d | 100 | 1,297 | – / 11 | n/a | 1,580.46* |
| X-n101-k25 (solved without MaxNV) | 100 | 206 | – / 25 | n/a | 27,591* |

*: Indicates the proven optimal objective value (*Opt*).

**Table A.2.** Temporal attributes of four time-restricted instances in the test bed.

| Instance Name | N | $D_{max}$ | fst |
|---|---|---|---|
| CMT-p06 | 50 | 200 | 10 |
| CMT-p07 | 75 | 160 | 10 |
| CMT-p08 | 100 | 230 | 10 |
| CMT-p14 | 100 | 1,040 | 90 |

**Table A.3.** Update for five *BKS* values posted on the VRP Web.

| Instance Name | True *Opt* | Incorrect or out-of-date *BKS* on VRP Web |
|---|---|---|
| ulysses-n16-k3 | 8,232 | 30,492 |
| P-n22-k8 | 603 | 590 |
| ulysses-n22-k4 | 9,312 | 40,153 |
| P-n55-k15 | 989 | 945 |
| E-n101-k14 | 1,067 | 1071 |

# References

Baldacci, R., Hadjiconstantinou, E., Mingozzi, A., (2004). An exact algorithm for the capacitated vehicle routing problem based on a two-commodity network flow formulation. Operations Research **52**(5), 723−738.

CVRPLIB (2015). Capacitated Vehicle Routing Problem Library. Algorithms, Optimization and Simulation Group (GALGOS), Pontifical Catholic University of Rio de Janeiro, Brasil. <http://vrp.atd-lab.inf.puc-rio.br/index.php/en/> *(accessed February 2016)*.

Dantzig, G., Fulkerson, R., Johnson, S. (1954). Solution of a large-scale traveling-salesman problem. Journal of the Operations Research Society of America **2**(4), 393−410.

Dantzig, G. B., Ramser, J. H. (1959). The truck dispatching problem. Management Science **6**(1), 80−91.

Desrochers, M., Laporte, G. (1991). Improvements and extensions to the Miller-Tucker-Zemlin subtour elimination constraints. Operations Research Letters **10**(1), 27−36.

Finke, G., Claus, A., Gunn, E. (1984). A two-commodity network flow approach to the traveling salesman problem. Congressus Numerantium **41**(1), 167−178.

GAMS Solver Manuals (November 2013). GAMS Version 24.2.1. GAMS Development Corporation. Washington, DC, USA.

Gavish, B., Graves, S. C. (1978). The travelling salesman problem and related problems. In: Working Paper GR-078-78, Operations Research Center, Massachusetts Institute of Technology.

Intel, 2015. Xeon® Processor W3690 Specifications. <http://ark.intel.com/products/52586/Intel-Xeon-Processor-W3690-12M-Cache-3_46-GHz-6_40-GTs-Intel-QPI> *(accessed February 2016)*.

Lenstra, J. K., Kan, A. H. G. (1981). Complexity of vehicle routing and scheduling problems. Networks **11**(2), 221−227.

Letchford, A. N., Salazar-González, J. J. (2006). Projection results for vehicle routing. Mathematical Programming, **105**(2-3), 251-274.

Letchford, A. N., Salazar-González, J. J. (2015). Stronger multi-commodity flow formulations of the Capacitated Vehicle Routing Problem. European Journal of Operational Research, **244**(3), 730-738.

Miller, C. E., Tucker, A. W., Zemin, R. A. (1960). Integer programming formulation of traveling salesman problems. Journal of the Association for Computing Machinery **7**(4), 326−329.

Öncan, T., Altınel, İ. K., Laporte, G. (2009). A comparative analysis of several asymmetric traveling salesman problem formulations. Computers & Operations Research **36**(3), 637−654.

Ordóñez, F., Sungur, I., Dessouky, M. (2007). A priori performance measures for arc-based formulations of vehicle routing problem. Transportation Research Record: Journal of the Transportation Research Board **2032**, 53−62.

Padberg, M., Sung, T. Y. (1991). An analytical comparison of different formulations of the travelling salesman problem. Mathematical Programming **52**(1−3), 315−357.

Roberti, R., Toth, P. (2012). Models and algorithms for the Asymmetric Traveling Salesman Problem: an experimental comparison. EURO Journal on Transportation and Logistics **1**(1–2), 113–133.

Sarin, S. C., Sherali, H. D., Judd, J. D., Tsai, P. F. J. (2014). Multiple asymmetric traveling salesmen problem with and without precedence constraints: Performance comparison of alternative formulations. Computers & Operations Research **51**, 64–89.

Toth, P., Vigo, D. (2002). The Vehicle Routing Problem. In: volume 9 of SIAM Monographs on Discrete Mathematics, SIAM Philadelphia.

Toth, P., Vigo, D. (2002). Models, relaxations and exact approaches for the capacitated vehicle routing problem. Discrete Applied Mathematics **123**(1), 487–512.

Uchoa, E., Pecin, D., Pessoa, A., Poggi, M., Subramanian, A., Vidal, T. (2014). New benchmark instances for the capacitated vehicle routing problem. Report downloaded from Optimization-Online.org, an eprint site for the optimization community. `<http://www.optimization-online.org/DB_HTML/2014/10/4597.html>` *(accessed February 2016)*.

VRP Web (2015). Capacitated VRP Instances. Networking and Emerging Optimization Research Group (NEO), University of Malaga, Spain. `<http://neo.lcc.uma.es/vrp/vrp-instances/capacitated-vrp-instances/>` *(accessed February 2016)*.
32